\documentclass{article}
\usepackage[utf8]{inputenc}

 \usepackage{tabu}

\usepackage{amsmath}

\usepackage{algorithm2e}
\usepackage{xcolor}

\SetAlFnt{\footnotesize}

\SetAlCapSty{xAlCapSty}

\SetCommentSty{xCommentSty}

\SetNlSty{mynlfont}{}{} 

\LinesNumbered

\SetSideCommentRight

\DontPrintSemicolon

\RestyleAlgo{algoruled}

\usepackage[driverfallback=dvipdfm,
bookmarks=true,    
bookmarksopen=true,
bookmarksopenlevel=\maxdimen,
bookmarksdepth=10
]{hyperref}
\hypersetup{
pdfcenterwindow=true, 
pdfpagelabels=true,
pagebackref=true,            
hypertexnames=true,
plainpages=false,
unicode=false,     
pdftoolbar=true,                
    pdfmenubar=true,          
    pdffitwindow=false,         
    pdfstartview={FitH},        
    pdftitle={On Third-Order Accuracy of MUSCL Scheme},
    pdfauthor={Hiroaki Nishikawa},    
    pdfsubject={On Third-Order Accuracy of MUSCL Scheme},       
    pdfcreator={Hiroaki Nishikawa},   
    pdfproducer={Hiroaki Nishikawa}, 
    baseurl={http://www.hiroakinishikawa.com},
    pdfkeywords={}, 
    pdfnewwindow=true,     
    colorlinks=true,       
    linkcolor=black,       
    citecolor=black,       
    filecolor=black,       
    urlcolor=black,        
  breaklinks=true,
  hyperfigures=true,
  backref=true,
breaklinks=false, 
pdfpagelabels,      
pagebackref,        
hypertexnames=true, 
plainpages=false,   
naturalnames        
}
\usepackage[all]{hypcap}

\usepackage[thinlines]{easytable}

\usepackage{color}      
\usepackage{placeins}
\usepackage{indentfirst}
\usepackage{multirow}

\usepackage{amssymb}
\usepackage{amsmath,latexsym}
\usepackage{url}
\usepackage{enumerate}
\usepackage{graphicx}

\usepackage{booktabs}

\usepackage{subcaption}
\usepackage{cleveref}
\usepackage{mathtools}

\usepackage{blindtext,graphicx}
\usepackage[absolute]{textpos}
\begin{document}


\title{\bf  Arithmetic Averages of Viscous Coefficients are Sufficient for Second-Order Finite-Volume Viscous Discretization on Unstructured Grids}

\author{ 
Hiroaki Nishikawa and Boris Diskin\\
  {\normalsize\itshape National Institute of Aerospace, Hampton, VA 23666, USA}
}

\date{\today}
\maketitle

\begin{abstract}  
In this short note, we discuss the use of arithmetic averages for the evaluation of viscous coefficients such as temperature and velocity components at a face as required in a cell-centered finite-volume viscous discretization on unstructured grids, and show that second-order accuracy can be achieved even when the arithmetic average is not linearly-exact second-order reconstruction at a face center (e.g., the face center is not located exactly halfway between two adjacent cell centroids) as typical in unstructured grids. Unlike inviscid discretizations, where the solution has to be reconstructed in a linearly exact manner to the face center for second-order accuracy, the viscous discretization does not require the linear exactness for computing viscous coefficients at a face. There are two requirements for second-order accuracy, and the arithmetic average satisfies both of them. Second-order accuracy is numerically demonstrated for a simple one-dimensional nonlinear diffusion problem and for a three-dimensional viscous problem based on methods of manufactured solutions. 
\end{abstract}


\section{Introduction}

We focus on unstructured tetrahedral grids and discretize the steady compressible Navier-Stokes (NS) equations
 by a second-order cell-centered finite-volume discretization, where the residual at a cell $j$ is given by
\begin{eqnarray}
{\bf Res}_j =  \sum_{k \in \{ k_j \}} {\bf \Phi}_{jk} |{\bf n}_{jk}| ,
\label{res_used}
\end{eqnarray}
where $\{ k_j \}$ is a set of neighbor cells, $ |{\bf n}_{jk}|$ is the area of the face between the cell $j$ and a neighbor $k$, ${\bf n}_{jk}$ is the scaled outward face-normal vector, $\Phi_{jk} $ is a numerical flux defined by ${\bf \Phi}_{jk} = {\bf \Phi}_{jk}^{inv} + {\bf \Phi}_{jk}^{vis}$, approximating the physical NS flux projected along $\hat{\bf n}_{jk} = {\bf n}_{jk} /|{\bf n}_{jk}|$, $ {\bf \Phi}_{jk}^{inv} $ and $ {\bf \Phi}_{jk}^{vis} $ denote the inviscid and viscous numerical fluxes, respectively. The numerical fluxes are computed at the centroid ${\bf x}_c$ of each triangular face using the primitive variables ${\bf w} = ( \rho,  {\bf v} , T )$, where $\rho$ is the density, $ {\bf v}  = (u,v,w)$ is the velocity vector, and $T$ is the temperature. The inviscid flux is computed using a Riemann solver and the left and right states, $ {\bf w}_L$ and $ {\bf w}_R$, linearly reconstructed from the cell centroids ${\bf x}_j$ and ${\bf x}_k$, respectively: e.g., by the Roe flux \cite{Roe_JCP_1981}, 
\begin{eqnarray}
  {\bf \Phi}_{jk}^{inv}  =    {\bf \Phi}_{jk}^{Roe}   ( {\bf w}_L, {\bf w}_R),
\label{numerical_flux}
\end{eqnarray}
where  
\begin{eqnarray}
  {\bf w}_L =   {\bf w}_j + \nabla {\bf w}_j \cdot  ( {\bf x}_c - {\bf x}_j ), \quad
  {\bf w}_R =   {\bf w}_k + \nabla {\bf w}_k \cdot ( {\bf x}_c - {\bf x}_k ),
\label{states_LR}
\end{eqnarray}
and the gradients, $\nabla {\bf w}_j $ and $\nabla {\bf w}_k$, are computed by a linear least-squares method. The resulting inviscid discretization is second-order accurate because the one-point flux quadrature is linearly exact over a triangular face. Here, our focus is on the viscous discretization represented by a general framework of evaluating the physical viscous flux with a face gradient $\nabla {\bf w}_f$ and face values $ {\bf v}_f$ and ${T}_f$: 
\begin{eqnarray} 
 {\bf \Phi}_{jk}^{vis}  = {\bf f}_n^{vis} ( {\bf w}_j, {\bf w}_k, \nabla {\bf w}_f ),   \quad
   {\bf f}_n^{vis}  = ( 0, - {\boldsymbol \tau}_{\!\! n},  - {\boldsymbol \tau}_{\!\!n}  {\bf v}_f +  {q}_n  ),
 \end{eqnarray} 
 where $ {\bf f}_n^{vis} $ is the physical viscous flux projected along $\hat{\bf n}_{jk}$ with zero, $ - {\boldsymbol \tau}_{\!\! n}$ , $- {\boldsymbol \tau}_{\!\!n}  {\bf v}_f +  {q}_n$ for the continuity, momentum, and energy equations, respectively, and 
\begin{eqnarray} 
q_n =  - \frac{   {\mu}_f   }{  P_r (\gamma - 1)}    \nabla \,T    \cdot  \hat{\bf n}_{jk},  
\quad
 {\boldsymbol \tau}_{\!\! n} = - {\mu}_f   \left[  \frac{2}{3} tr(  \nabla {\bf v} ) {\bf I} +   \nabla {\bf v}  + (  \nabla {\bf v})^t  \right]  \hat{\bf n}_{jk} ,
 \quad
 {\mu}_f  =   \frac{M_{\infty}}{Re_{\infty}}  \frac{1+C /{T}_{\infty}}{ {T}_f +C/ {T}_{\infty}}  {T}_{\! f}^{\frac{3}{2}} .
\end{eqnarray} 
Here, $tr()$ denotes the trace, the superscript $t$ indicates the transpose, ${\bf I}$ is the 3$\times$3 identity matrix, $M_{\infty}$ is a free stream Mach number, $Re_{\infty}$ is a free stream Reynolds number, ${T}_{\infty}$ is a dimensional free stream temperature, and $C = 110.5 $ [K] is the Sutherland constant. All the quantities are assumed to have been nondimensionalized by their free-stream values except that the velocity and the pressure are scaled by the free-stream speed of sound and the free-stream dynamic pressure, respectively, which has led to $p = \rho T / \gamma$ and the factor $M_{\infty} / Re_{\infty}$ in the viscosity (see Ref.~\cite{idolikeCFD_VOL1_2013_pdf}). 
The choice of the face gradient is irrelevant to the discussion, but for the numerical results in this short note, we employed the alpha-damped face gradient formula \cite{nishikawa:AIAA2010}: 
\begin{eqnarray}
 \nabla {\bf w}  =   \frac{1}{2} \left[   \nabla^{LSQ} {\bf w}_j + \nabla^{LSQ}  {\bf w}_k  \right] + \frac{\alpha}{ | ( {\bf x}_k -   {\bf x}_j ) \cdot \hat{\bf n}_{jk}| } 
({\bf w}_R -{\bf w}_L) \hat{\bf n}_{jk},  \quad  
\alpha=4/3.
\end{eqnarray}

To complete the viscous disretization, we need to define the face quantities: $ {T}_{\! f}$ and ${\bf v}_f$; this is the main focus of this short note. Since reconstructed solutions are available at the face center, we may evaluate them as
\begin{eqnarray} 
 {T}_{\! f} = \frac{ T_L  + T_R }{2}, \quad {\bf v}_f = \frac{  {\bf v}_L + {\bf v}_R }{2}.
 \label{Tuvw_aveLR}
\end{eqnarray}
These formulas are linearly exact on arbitrary tetrahedral grids and achieve second-order accuracy. However, in many practical cell-centered computational fluid dynamics (CFD) codes (e.g., those presented in Refs. \cite{VULCAN_FANGplus_scitech2022,pandya_etal:AIAAJ2016,scFLOW:Aviation2020} although not explicitly mentioned), these quantities are evaluated with the arithmetic average of the cell values:
\begin{eqnarray} 
 {T}_{\! f}  = \frac{ T_j  + T_k }{2}, \quad  {\bf v}_f = \frac{ {\bf v}_j +  {\bf v}_k }{2}.
 \label{Tuvw_avejk}
\end{eqnarray}
This reconstruction is not linearly exact if the face center is not located halfway between the two adjacent cell centroids, but it can be more robust since the face temperature is guaranteed to be positive as long as the cell values of the temperature are positive, for example. The point of this short note is that the arithmetic average is actually sufficient for achieving second-order accuracy in the viscous discretization even if the reconstruction is not linearly exact.

\section{Two Requirements for Second-Order Accuracy}

The first and primary requirement for second-order accuracy in a finite-volume viscous discretization is that the gradients be computed by an algorithm that is exact for linear functions on irregular grids  \cite{Boris_Jim_NIA2007-08}, or in other words, that the viscous flux at a face be first-order accurate on irregular grids. This requirement is typically met by the use of a linear least-squares method for computing gradients at cells. The reconstruction of a face quantity used in the viscosity coefficients however, does not have to be linearly exact but just have to be first-order accurate because first-order errors are already committed in the gradients that the viscosity coefficients multiply. For example, for the arithmetic average, we find
\begin{eqnarray} 
  \frac{ T_j  + T_k }{2}  
  &=&   \frac{ T({\bf x}_c) + ({\bf x}_c-{\bf x}_j) \cdot \nabla T({\bf x}_c)  +  T({\bf x}_c) + ({\bf x}_c - {\bf x}_k) \cdot \nabla T({\bf x}_c) }{2}  +O(h^2) 
\nonumber  \\ [2ex]
  &=& T({\bf x}_c) +  \left(  {\bf x}_c - \frac{ {\bf x}_j +  {\bf x}_k}{2}  \right) \cdot \nabla T({\bf x}_c)  +O(h^2) 
\nonumber  \\ [2ex]
  &=& T({\bf x}_c) + O(h),
\end{eqnarray}
on unstructured grids where ${\bf x}_c \ne \frac{ {\bf x}_j +  {\bf x}_k}{2}$. Then, the viscous stresses are the product of the face viscosity evaluated with the face temperature and the velocity gradients, and thus they are first-order accurate overall. Therefore, a first-order accurate reconstruction is sufficient (if not linearly exact) for achieving second-order accuracy in the viscous discretization on general unstructured tetrahedral grids. The same is true for other types of grids if applied with a linearly exact flux quadrature (e.g., split a quadrilateral face into two triangles and apply the one-point quadrature). The argument equally applies to the face velocity ${\bf v}_f$ at a face required in the energy equation. 

For practical purposes, the first-order accurate reconstruction is necessary, but not sufficient. The following one-sided evaluation
\begin{eqnarray} 
{T}_f  = T_j, 
\end{eqnarray}
 is also first-order accurate but does not lead to second-order accuracy. The additional requirement is that the reconstruction has to become linearly exact on regular grids. This requirement excludes reconstructions that are biased to one side. In the next section, we will present numerical results to demonstrate second-order accuracy with the arithmetic average and first-order accuracy with one-sided evaluation reconstructions.

\section{Results}
\label{sec:results}

\subsection{Nonlinear Diffusion in One Dimension}
\label{sec:results_oned_nonlinear_diff} 

\begin{figure}[htbp!]
    \centering
      \begin{subfigure}[t]{0.32\textwidth}
        \includegraphics[width=\textwidth,trim=0 0 0 0 clip]{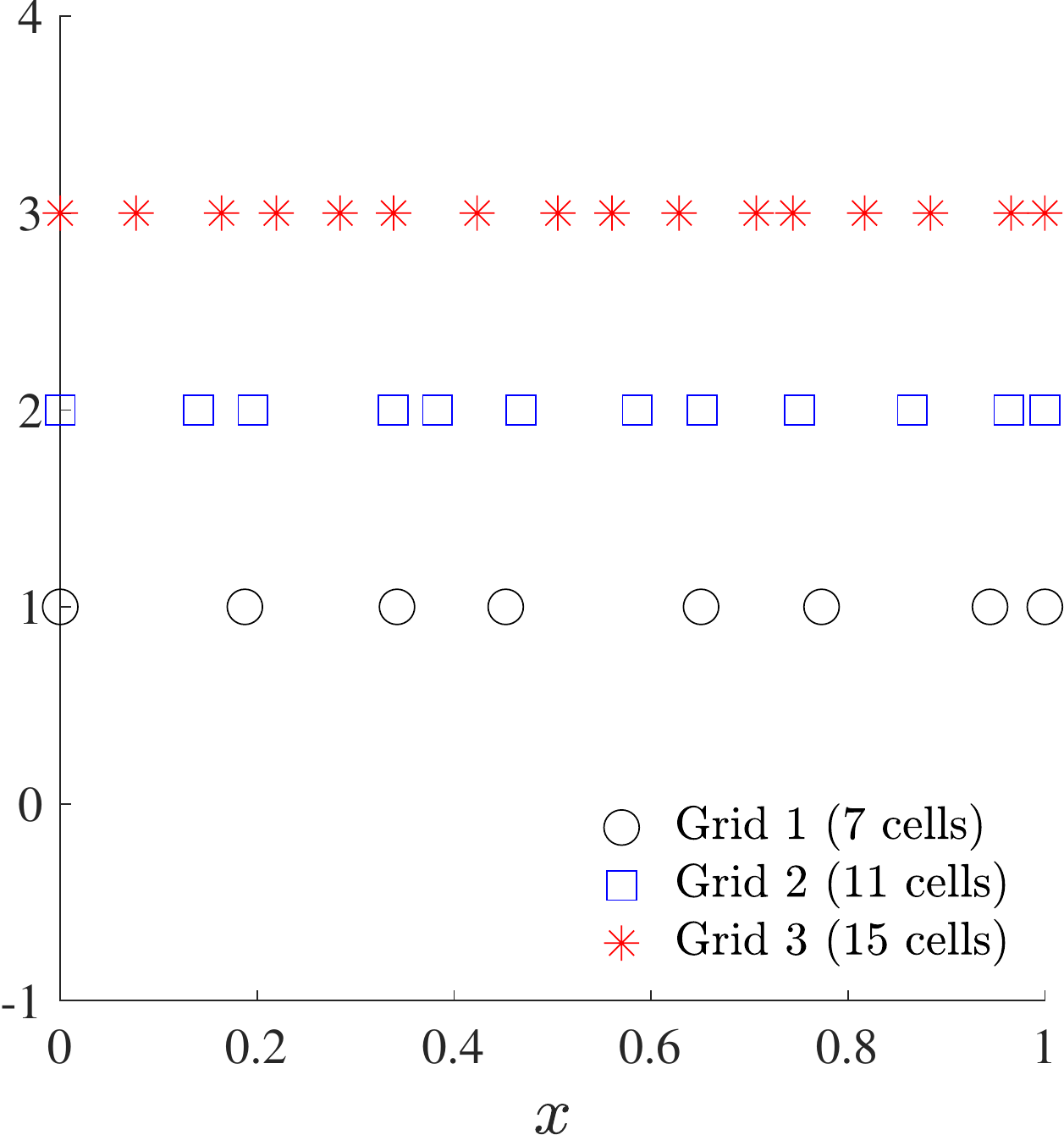}
          \caption{Three grids from the coarsest one.}
          \label{fig:oned_grid}
      \end{subfigure}
      \hfill
      \begin{subfigure}[t]{0.32\textwidth}
        \includegraphics[width=\textwidth,trim=0 0 0 0 ,clip]{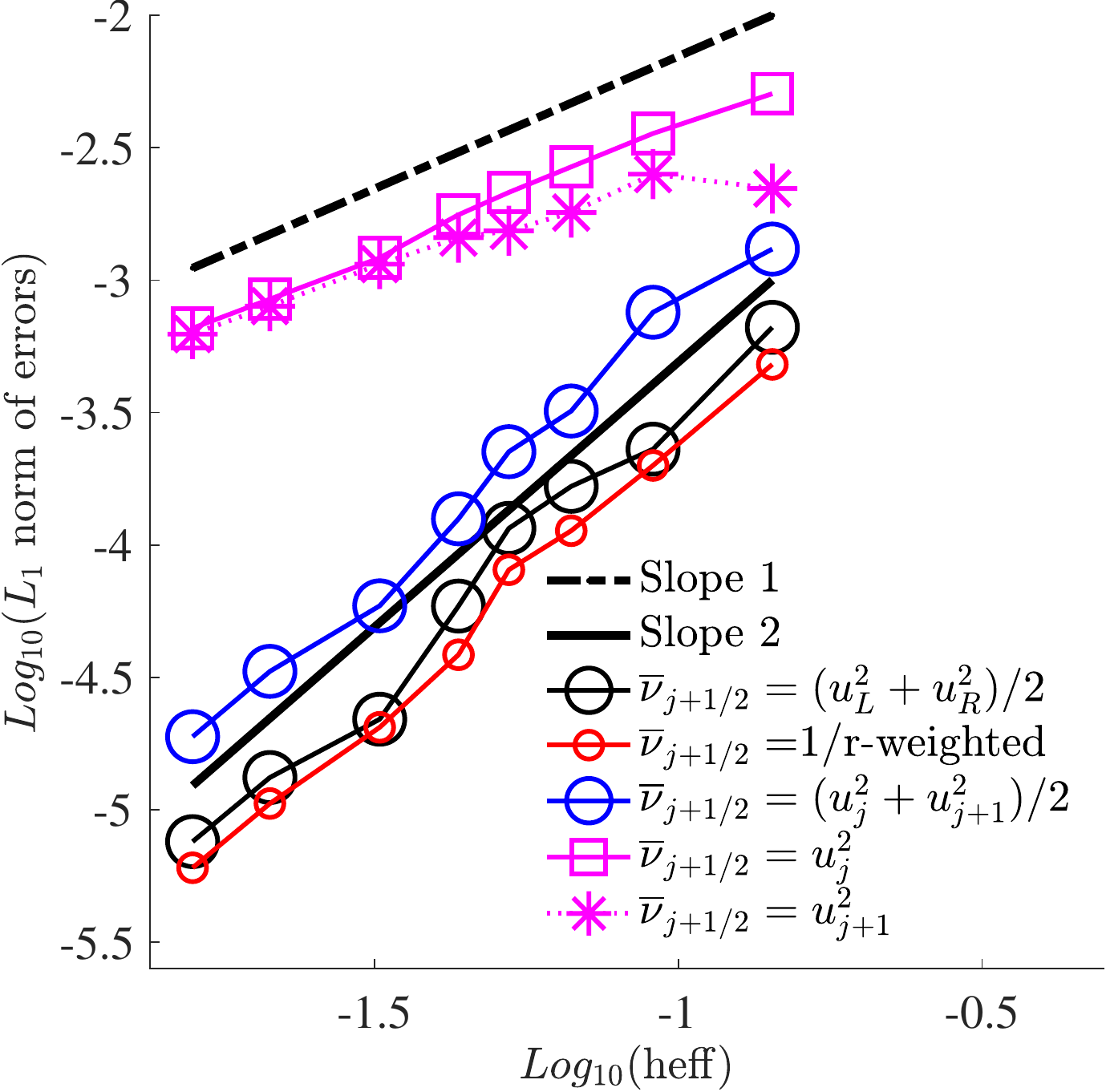}
          \caption{Irregular grids.}
          \label{fig:oned_de_irreg}
      \end{subfigure}
      \hfill
      \begin{subfigure}[t]{0.32\textwidth}
        \includegraphics[width=\textwidth,trim=0 0 0 0 ,clip]{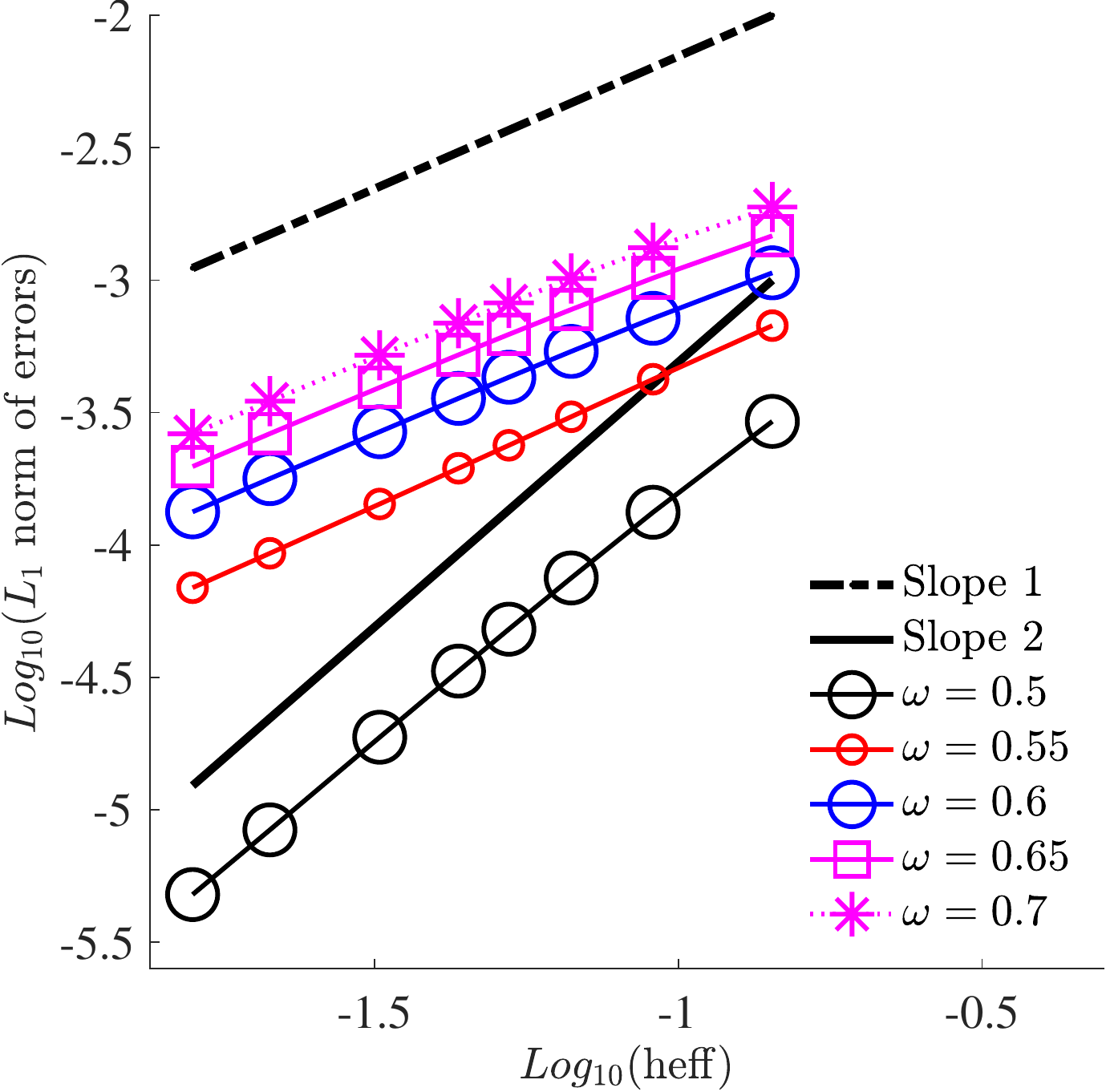}
          \caption{Regular grids.}
          \label{fig:oned_de_reg}
      \end{subfigure}
            \caption{ Grids and convergence of discretization errors for the one-dimensional problem: (a) symbols indicate node/face locations (a cell center is defined as the midpoint of two consecutive nodes), (b) error convergence for five different viscosity evaluations, (c) error convergence for a weighted-averaged viscosity with varying weight $\omega$.  }
\label{fig:oned}
\end{figure}

We consider a nonlinear diffusion problem in $x \in [0,1]$:
\begin{eqnarray}
  -  \partial_x \left( \nu \partial_x u  \right) = f(x), \quad \nu = u^2, \quad f(x) =   -  \partial_x \left( \nu \partial_x u_e  \right) ,
\end{eqnarray}
where the exact solution $u_e$ is given by
\begin{eqnarray}
  u_e  = \exp ( 2 x ).
\end{eqnarray}
We solve this problem on eight levels of irregularly-spaced grids with $n=7$, $11$, $15$, $19$, $23$, $31$, $47$, and $63$ cells (see Figure \ref{fig:oned_grid}) by a cell-centered finite-volume method. The discrete residuals are given, for $j=1,2,3, \cdots, n$, 
\begin{eqnarray}
 Res_j = \phi_{j+1/2} - \phi_{j-1/2} - f(x_j) h_j,
\end{eqnarray}
where $x_j$ denotes the cell center coordinate of the cell $j$, $h_j$ is the cell volume, and the numerical flux is given by
\begin{eqnarray}
   \phi_{j+1/2}  =  {\nu}_{j+1/2}  (u_x)_{j+1/2},    \quad
    (u_x)_{j+1/2} =    \frac{1}{2}   [   (u_x)_j + (u_x)_{j+1}  ]  + \frac{\alpha}{ 2 ( x_{j+1} - x_j  )}  (u_R - u_L) , \\ [1.5ex]
    u_L = u_j +  (u_x)_j  ( x_f - x_j ), \quad  u_R = u_{j+1} +  (u_x)_{j+1}  ( x_f - x_{j+1}  ) , \\ [1.5ex]
     (u_x)_j   = \frac{u_{j+1} - u_{j-1} }{  x_{j+1} - x_{j-1}  } ,  \quad    (u_x)_{j+1} = \frac{u_{j+2} - u_{j} }{  x_{j+2} - x_j  },   
\end{eqnarray}
with $\alpha=4/3$ \cite{nishikawa:AIAA2010}, where $x_f = x_{j + 1/2}$ denotes the face location. For simplicity, we specify the exact solution in the cells adjacent to the boundaries: 
\begin{eqnarray}
u_1 = u_e(x_1), \quad u_n = u_e(x_n), \quad Res_1=0, \quad Res_n=0. 
\end{eqnarray}
The system of discrete equations is solved by an implicit defect-correction solver until the residual is reduced by eight orders of magnitude from the initial value in the $L_1$ norm.

To investigate the impact of the reconstruction of the face quantities on solution accuracy, we evaluate the viscosity at a face in five different ways:
\begin{eqnarray}
    {\nu}_{j+1/2}  = \frac{ u_L^2 + u_R^2 }{2}, \quad
    {\nu}_{j+1/2}  = \frac{   ( x_f - x_j )^{-1} u_j^2 +  ( x_f - x_{j+1} )^{-1}  u_k^2 }{  ( x_f - x_j )^{-1}  +( x_f - x_{j+1} )^{-1} },  \label{oned_vis_01}\\ [1ex]
    {\nu}_{j+1/2}  = \frac{ u_j^2 + u_{j+1}^2 }{2}, \quad
    {\nu}_{j+1/2}  = u_j^2 , \quad
    {\nu}_{j+1/2}  = u_{j+1}^2.
\end{eqnarray}
Grid convergence plots of discretization errors are shown in Figure \ref{fig:oned_de_irreg}. The inverse-distance weighted reconstruction as in the second equation in Equation (\ref{oned_vis_01}) is referred to as $1/r$-weighted. As expected, second-order accuracy is achieved with $  {\nu}_{j+1/2}  = \frac{ u_L^2 + u_R^2 }{2}$, the inverse-distance weighting, and the arithmetic average  $ {\nu}_{j+1/2}  = \frac{ u_j^2 + u_{j+1}^2 }{2}$. On the other hand, the accuracy order is deteriorated to first-order with the two one-sided viscosity evaluations: ${\nu}_{j+1/2}  = u_j^2$ and ${\nu}_{j+1/2}  = u_{j+1}^2$. To illustrate this accuracy deterioration problem, we performed computations with regular grids of the same sizes and with the following weighted average: 
\begin{eqnarray}
    {\nu}_{j+1/2}  =   \omega u_j^2 + (1-\omega)  u_{j+1} ^2,
\end{eqnarray}
which reduces to ${\nu}_{j+1/2}  = \frac{ u_j^2 + u_{j+1}^2 }{2}$ at $\omega=1/2$ and to ${\nu}_{j+1/2}  = u_j^2$ at $\omega=1$. 
The grid convergence plots  are shown in Figure \ref{fig:oned_de_reg}. Second-order accuracy is achieved only with $\omega=1/2$, and all others lead to first-order accuracy, which was expected because the weighted average is linearly exact on these regular grids only at $\omega=1/2$. These results support the two requirements discussed earlier: (1) a diffusive flux is at least first-order accurate and (2) the face reconstruction is linearly exact on regular grids.

%

\subsection{Viscous Problem in Three Dimensions}
\label{sec:results_visc_3d} 

\begin{figure}[htbp!]
    \centering
      \begin{subfigure}[t]{0.48\textwidth}
        \includegraphics[width=\textwidth,trim=3 3 3 3,clip]{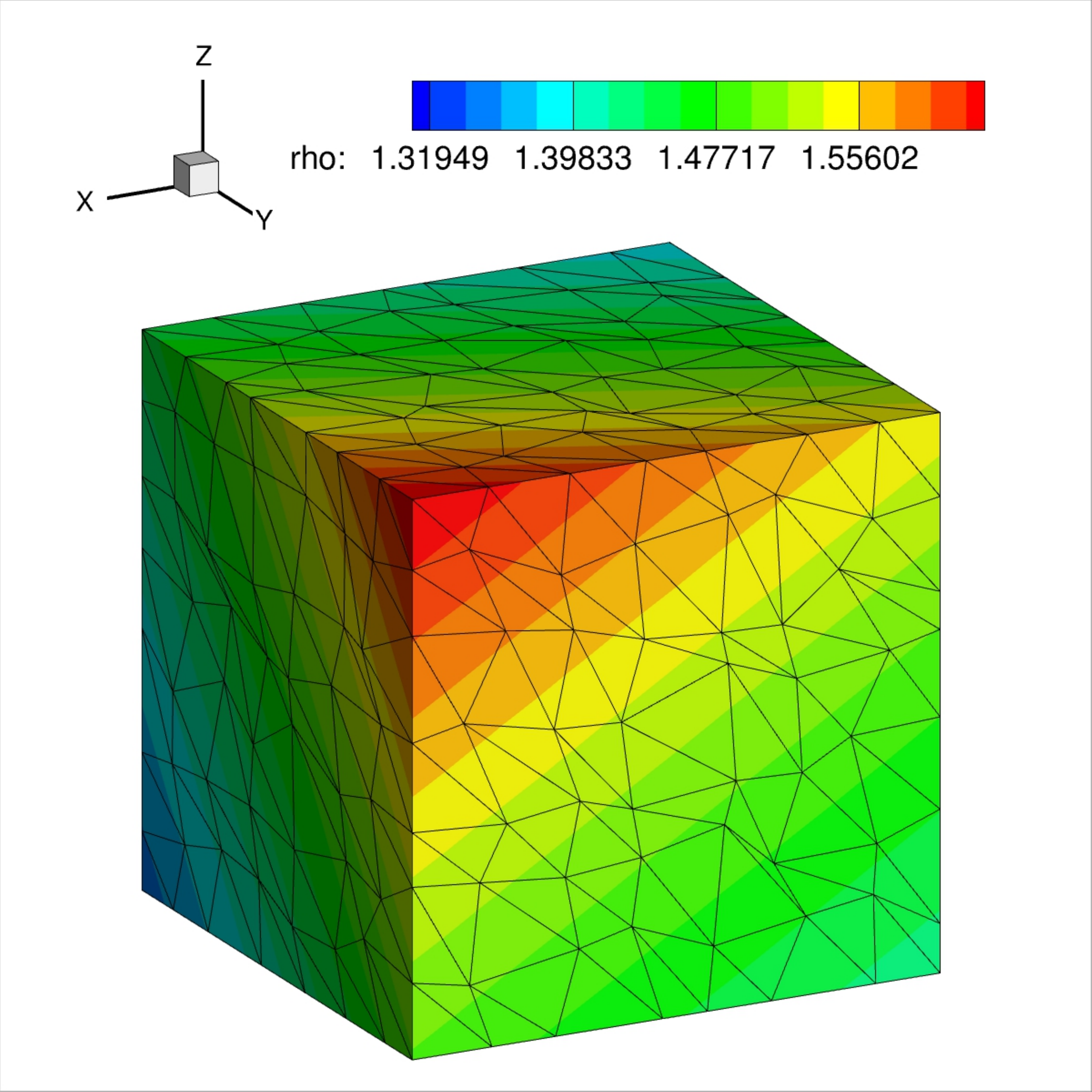}
          \caption{Coarsest grid.}
          \label{fig:3d_grid}
      \end{subfigure}
      \hfill
      \begin{subfigure}[t]{0.48\textwidth}
        \includegraphics[width=\textwidth,trim=0 0 0 0 clip]{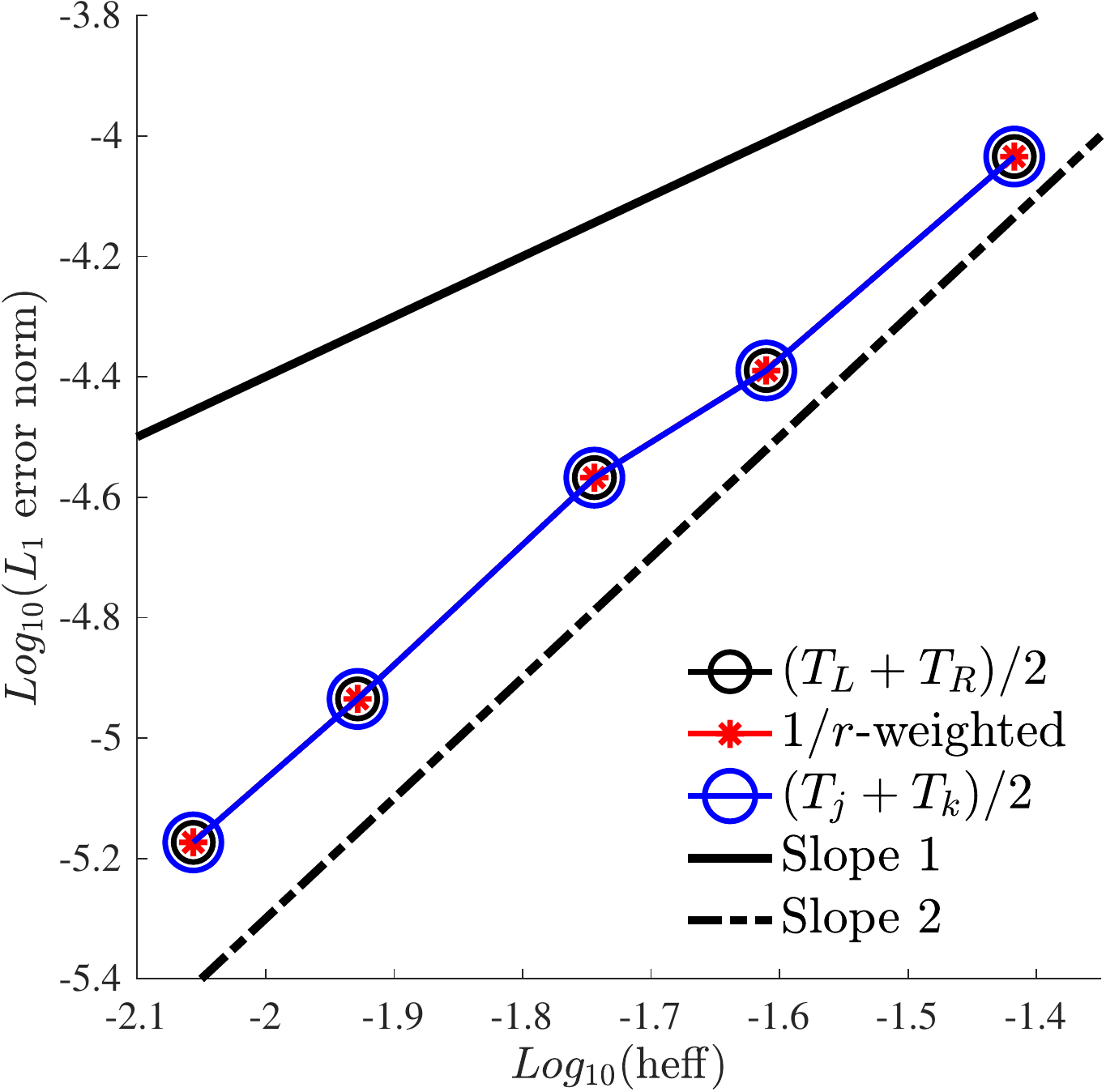}
          \caption{Error convergence for $\rho$.}
          \label{fig:3d_de_irreg}
      \end{subfigure} 
            \caption{ Grids and results for the three-dimensional problem: (a) the coarsest grid (all the boundaries are flat), (b) error convergence results for the density.  }
\label{fig:3d}
\end{figure}

In this section, we verify second-order accuracy with the arithmetic average and inverse-distance reconstructions for the compressible NS system in three dimensions by using the method of manufactured solutions. We introduce a vector of forcing terms such that the exact solution is given by
\begin{eqnarray}
{\bf w} = ( \rho,  {\bf v} , T ) = 
\left(
  1.0 ,
  0.3 ,
  0.2  ,
  0.1  ,
  1.0  
\right) 
+
\left(
  1,1,1,1,1 
  \right) \psi(x,y,z),
\end{eqnarray}  
where $ \psi(x,y,z) =  0.1 \exp( 0.5 (x + y + z) )$. The NS system is discretized in a cube domain ($(x,y,z) = [0,0.5] \times [0,0.5] \times [0,0.5] \times $) by the cell-centered finite-volume discretization with a point evaluation of the forcing term vector and three different face reconstruction methods: Equation (\ref{Tuvw_aveLR}), Equation (\ref{Tuvw_avejk}), and the inverse-distance weighted reconstruction,
\begin{eqnarray} 
   T_f  = \frac{   | {\bf x}_f - {\bf x}_j |^{-1} T_j+  | {\bf x}_f - {\bf x}_k |^{-1} T_k }{   | {\bf x}_f - {\bf x}_j |^{-1}  +  | {\bf x}_f - {\bf x}_k |^{-1} }. 
   \label{idw_3d}
\end{eqnarray}
 A family of six grids with 2058, 7986, 20250, 73002, and 178746 cells has been generated. Grids are all tetrahedral and irregular, and each grid has been generated independently. The coarsest grid is shown in Figure \ref{fig:3d_grid}. 
For simplicity, as in the one-dimensional case, we specify the exact solution in the cells adjacent to the boundary. To study the impact of the viscous discretization, we set $M_\infty$=0.1 and $Re_\infty = 0.1$ with $T_\infty =300 [K]$. The nonlinear discrete equations are solved by an implicit defect-correction solver until the residual is reduced by five orders of magnitude in the $L_1$ norm.

The discretization error convergence plots are shown in Figure \ref{fig:3d_de_irreg} for the density (results are similar for other primitive variables and therefore not shown). Clearly, all the face reconstruction methods achieve second-order accuracy. Actual errors are different but very similar: e.g., for the finest grid, the $L_1$ norm of the discretization error is $6.70945$e-$05$ for $T_f = (T_L + T_R)/2$, 
$6.70906$e-$06$ for $T_f = (T_j + T_k)/2$, and $6.70911$e-$06$ for the inverse-distance weighted reconstruction (\ref{idw_3d}).

%

\section{Conclusions}
\label{conclusions}
 
We have discussed the use of the arithmetic average for the face quantities in the viscous flux within a cell-centered finite-volume viscous discretization for unstructured grids, and have shown that the arithmetic average is sufficient for achieving second-order accuracy in the viscous discretization even if the face center is not located at halfway between two adjacent cells. The study suggests two requirements for the evaluation of face quantities for the viscous flux: (1) the viscous flux is first-order accurate, and (2) the face reconstruction is linearly exact on regular grids. Numerical results have been presented to illustrate the importance of these requirements. 
 
This short note is prepared for providing an explanation for second-order accuracy that is achieved with the arithmetically-averaged face quantities in the viscous flux on unstructured grids. Many practical unstructured-grid CFD codes have already been using the arithmetic average \cite{VULCAN_FANGplus_scitech2022,pandya_etal:AIAAJ2016,scFLOW:Aviation2020}, and it will continue to be used with confidence as it maintains second-order accuracy. 

\addcontentsline{toc}{section}{Acknowledgments}
\section*{Acknowledgments}

The first author gratefully acknowledges support by the Hypersonic Technology Project, through the Hypersonic Airbreathing Propulsion Branch of the NASA Langley Research Center, under Contract No. 80LARC17C0004.
 

\addcontentsline{toc}{section}{References}
\bibliography{../../bibtex_nishikawa_database}
\bibliographystyle{aiaa}


\end{document}